\documentclass[11pt]{article}

\usepackage{amsmath, amssymb}
\usepackage{graphicx}
\usepackage{float}

\newcommand{\cqfd}{{\nobreak\hfil\penalty50\hskip2em\hbox{}
\nobreak\hfil $\square$\qquad\parfillskip=0pt\finalhyphendemerits=0\par\medskip}}
\newcommand{\R}{\mathbb{R}}
\newcommand{\N}{\mathbb{N}}

\newcommand{\dt}{\partial_t}
\newcommand{\dx}{\partial_x}
\newcommand{\ds}{\displaystyle}

\newcommand{\pr}{{\bf \textit{Proof : }}}

\def\bR{\mathbf{R}}
\def\bS{\mathbf{S}}

\newtheorem{theorem}{Theorem}[section]

\newtheorem{proposition}{Proposition}[section]
\newtheorem{lemma}{Lemma}[section]
\newtheorem{definition}{Definition}[section]
\newtheorem{remark}{Remark}[section]

\title{ \bf 
Blow up at the hyperbolic boundary for a $2\times2$  system arising from chemical engineering}

\author {C. Bourdarias
\thanks{Universit\'{e} de Savoie, LAMA, UMR CNRS 5127, 73376 Le Bourget-du-Lac,
bourdarias@univ-savoie.fr},
M. Gisclon
\thanks{Universit\'{e} de Savoie, LAMA, UMR CNRS 5127, 73376 Le
Bourget-du-Lac, gisclon@univ-savoie.fr}
and S. Junca
\thanks{ Universit\'{e} de Nice, Labo. JAD, UMR CNRS 6621, Parc Valrose, 06108
Nice, junca@unice.fr}
}
\date{November 22, 2009}

\begin{document}

\bibliographystyle{plain}

\maketitle

\abstract{
We consider an initial boundary value problem for a 2x2  system of conservation laws mo\-deling
heatless adsorption of a gaseous mixture with two species and instantaneous exchange kinetics, close
to the system of Chromatography. In this model the velocity is not constant because the sorption
effect is taken into account. Exchanging the roles of the $x$, $t$ variables we obtain a strictly
hyperbolic system with a zero eigenvalue. Our aim is to construct a solution with a velocity which
blows up  at the corresponding characteristic ``hyperbolic boundary'' $\{t=0\}$.}
\medskip

\noindent {\bf AMS Classification}: 35L65, 35L67, 35Q35.

\medskip
 \noindent {\bf Key words}: 
hyperbolic systems, conservation laws, weak solutions, Temple systems, blow up, boundary conditions,
gas chromatography, Riemann problem, Front Tracking Algorithm.
%
\section{Introduction} 
%

This paper deals with the construction of a blow up  solution for a one dimensional 2x2 hyperbolic
system of conservation laws related to a particular isothermal gas-solid chromatography process,
called ``Pressure Swing Adsorption'' (PSA), with two species. Each of the two species simultaneously
exists under two phases: a gaseous and movable one with a common velocity $u(t,x)$ and concentration
$c_i (t,x) \geq 0$ and a solid (adsorbed)  with concentration $q_i$, $i=1, \, 2$. 
One may consult \cite {RAA70} and \cite{Ru84} for a precise description of the process and
\cite{BGJ08} for a survey on various related models. 

In gas chromatography, velocity variations accompany changes in gas composition, especially in the
case of high concentration solute: it is known as the sorption effect. Here, the sorption effect is
taken into account through a constraint on the pressure that is, in this isothermal model, on the
total density:
\begin{equation}\label{constraint}
 c_1+c_2=\rho(t),
\end{equation} 
where the function $\rho$ is given (this is actually achieved in the experimental or industrial
device). This constraint  expresses that the total pressure depends only on the time. In the sequel
we assume that $\rho\equiv 1$. 

We assume that mass exchanges between the mobile and the stationary phases are infinitely fast
(instantaneous exchange kinetics) thus the two phases are constantly at  composition
e\-qui\-li\-brium: the concentrations in the solid  phase are given by some relations 
$q_i=q_i^*(c_1,c_2)$ where the functions $q_i^* $ are the so-called e\-qui\-li\-brium isotherms.
A theoretical study  of a model with finite exchange kinetics was presented in \cite{B92} and a
numerical approach was developed in \cite{B98}. 

With these assumptions, the PSA system reads:
\begin{eqnarray}
\dt (c_1+q^*_1(c_1,c_2))+\dx(u\,c_1)&=&0, \label{un}\\
\dt (c_2+q^*_2(c_1,c_2))+\dx(u\,c_2)&=&0, \label{deux}\\
c_1+c_2&=& 1. \label{trois}
\end{eqnarray}
Notice that $c_1$, $c_2$ must satisfy $0\leq c_1,\,c_2\leq
1$. 

Adding (\ref{un}) and (\ref{deux}) we get, thanks to (\ref{trois}):
$$\dt (q^*_1(c_1,c_2)+q^*_2(c_1,c_2))+\dx u=0,$$

i.e. the constraint (\ref{constraint}) leads to an integral dependency of the velocity upon the
concentrations. We denote $c=c_1$ then $c_2=1-c$ and the unknowns are the velocity $u$ and the
concentration $c$. We write the PSA system under the form:
\begin{eqnarray}
\dx u  + \dt h(c) & =&  0, \label{psa1}\\
\dx(u\,c)+ \dt I(c) &  =&  0,\label{psa2}
\end{eqnarray}
with
\begin{eqnarray*}
h(c)&=&q^*_1(c,1-c)+q^*_2(c,1-c),\\
I(c)&=&c+q^*_1(c,1-c).
\end{eqnarray*}

In the sequel we denote $$q_1(c)=q^*_1(c,1-c)\quad\hbox{and}\quad q_2(c)=q^*_2(c,1-c),$$
thus $h=q_1+q_2$.
Any equilibrium isotherm related to a given species is always non decreasing with respect to the
corresponding concentration and non increasing with respect to the others (see
\cite{DCRBT88}) i.e. 
$$\frac{\partial q^*_i}{\partial c_i}\geq 0\hbox{ and }\frac{\partial q^*_i}{\partial c_j}\leq
0\hbox{ for }j\neq i,$$
thus we have immediately:
\begin{equation}\label{qprime}
q'_1\geq   0  \geq q'_2.
\end{equation}

PSA system (\ref{psa1})-(\ref{psa2}) is completed by initial and boundary data:
\begin{equation} \label{psa0}
\left\{   \begin{array}{ccl}
 \vspace{2mm}c(0,x)&=&c_0(x) \in [0,1],  \quad x > 0,\\
 \vspace{2mm}c(t,0) &=&c_b(t)   \in [0,1],\quad t>0,\\
 u(t,0)&=&u_b(t)\geq\alpha,\quad t>0,\quad\hbox{for some constant }\alpha>0.
\end{array}\right.
\end{equation}

Notice that we assume in (\ref{psa0}) an influx boundary condition, i.e. $\forall t>0,\,u_b(t)>0$
and we choose $]0,+\infty[$ instead of $]0,1[$ as spatial domain for the sake of simplicity.

It is well known that it  is possible to analyze the system of Chromatography, and thus PSA system,
in terms of hyperbolic system of P.D.E. provided to exchange the time and space variables: see
\cite{RAA86} and also \cite{RSVG88} for instance. In this framework,  with evolution variable $x$
instead of $t$, the system is strictly hyperbolic and has a null eigenvalue. Thus $\{t=0\}$ is a
characteristic boundary. 

We obtained in \cite{BGJ06,BGJ07} an existence result for weak global entropy solutions with large
$BV$ data for the concentrations and only $L^\infty$ data for the velocity: as in
\cite{CG99,CG01,M04}, the zero eigenvalue makes  the existence possible of stratified
solutions in the sense that $u(t,x)=u_b(t)\, v(t,x)$ where $v$ is as regular as the concentration
$c$ and more than  the ``initial'' data $u_b$ (see \cite{BGJ4}). 
This  linearly degenerate  eigenvalue  makes also possible the propagation  of large amplitude
and high frequency waves as in \cite{CGM04,CGM03}. 
As we will see below, we cannot expect in general to have weak global entropy solutions with
$L^\infty$ data for the concentrations because it is possible, in that case, to build a blow up
solution for some particular isotherms related for instance to  ammonia or water vapor.

It is already known that systems of two  hyperbolic conservation laws may blow up in the
$L^\infty$ norm. In literature (\cite{SY05,Y99}), there are examples of blow up for one dimensional
strictly hyperbolic systems of partial differential equations. For $2\times 2$ strictly hyperbolic
system there are few exam\-ples (\cite{S07}). In \cite{Y03short}, see also \cite{Y03},  R. Young  built  a very nice exam\-ple
involving two Burgers equations, linearly coupled at the two boundaries. All examples occur
in cases where strict hyperbolicity is lost as the solution explodes (see also  H.-K.Jenssen and  R.
Young \cite{JY04} page 630). In our exam\-ple, we also loose the strict hyperbolicity of the system at
the blow up point, but the blow up takes place only at the characteristic boundary $\{t=0\}$  which
becomes twice characteristic and only the velocity blows up. The context of our example (gas-solid
chromatography) is physically motivated and the various assumptions are physically relevant
(they are achieved for particular gases) but a blow up along the $x$-axis at $t=0$ is
of course artificial. Nevertheless it illustrates what may occur when $BV$ regularity is not
ensured for the velocity at the physical boundary.

The aim of this paper is to present this example of blow up for a class of
PSA systems. The paper is organised as follows. In Section \ref{secH}, we give an overview on
 PSA system (\ref{psa1})-(\ref{psa2}) and specify our  assumptions. In Section
\ref{secTemple} we investigate the conditions on PSA system (\ref{psa1})-(\ref{psa2}) to be or not
to be in the Temple class: indeed there is no blow up for such  systems. In Section \ref{secBlowup}
we state, prove and comment the main result of this paper about a blow up of velocity in $L^\infty$.

\section{The PSA system} \label{secH}

\subsection{General properties} 


In this section, we recall briefly some properties of  PSA system (\ref{psa1})-(\ref{psa2}). For
more details and explanations  we refer the reader to \cite{BGJ06,BGJ07,BGJ08,BGJ4}. 

As pointed out by Rouchon and \textit{al.} (\cite{RSVG88}), it is possible to analyse PSA system
in terms of hyperbolic system of P.D.E. provided we exchange the time and space variables. In this
framework the vector state will be  $U=\left(
\begin{array}{l}
u \\
m
\end{array}\right) $ where $m=u\,c $ is the flow rate of the first species. The first component $u$
of this vector must be understood as $u\,\rho$, that is the total flow rate. 
The reader should bear in mind that convexity is a relevant property of the entropy only with
respect to the conservative variables $(u,m)$ (see Dafermos \cite{D00} section 4.5 page 76  or 
\cite{BGJ07}). Notice that we seek for solutions satisfying $u>0$, wich ensure the system to be
hyperbolic. Using  the diffeomorphism $(u,c)\mapsto (u,m)$ between $\R_+^*\times ]0,1[$ and
$\R_+^*\times  \R_+^*$, it is more convenient  to deal with the hyperbolicity, the Riemann
invariants, the Riemann problem and the Front tracking algorithm with the  set of variables
$(u,c)$. For instance there is no  partial derivative with respect to $t$ for $u$ and $c$ is a
Riemann invariant.

PSA system  with conservative variables $(u,m)$  takes the form
\begin{equation}\label{sysadum}
\partial_x U +\partial_t \Phi(U)=0\hbox{  with }
\Phi(U)=\left(
\begin{array}{l}
h(m/u) \\\\
I(m/u)
\end{array}\right).
\end{equation}
It is strictly hyperbolic with variable $x$ as the evolution variable: the two
eigenvalues are $0$ (linearly degenerate) and $\lambda=H(c)/u$ where
\begin{equation}\label{eqH}
 H(c)=1+q'_1(c)-ch'(c)=1+(1-c)\,q'_1(c)-c\,q'_2(c)\geq 1.
\end{equation} 
In particular, the mathematical boundary condition  is  now related to the set
$\{t=0\}$ and the mathematical initial data to the set $\{x=0\}$.
 Making use of the function $f=q_1\,c_2-q_2\,c_1$ introduced by Douglas and {\em al.} in
\cite{DCRBT88}, written here under the form 
\begin{equation}  \label{eqf}
    f(c)= q_1(c)-c\,h(c),
\end{equation} 
it is easy to show that the right eigenvector 
$r=\left(\begin{array}{c} h'(c) \\
1+q'_1(c)
\end{array} \right )$ 
associated to $\lambda$ satisfies $\ds d\lambda \cdot r=\ds\frac{H(c)}{u^2}\,f''(c)$, thus $\lambda$
is genuinely nonlinear in each domain where $f''\neq 0$. 
In view of the remark following (\ref{trois}) there is a ``natural maximum principle'' for the
concentration $c$ : 
\begin{equation}\label{maxp}
 0  \leq c \leq 1.
\end{equation}
and this is achieved by the solutions (see \cite{BGJ06}).

PSA system (\ref{psa1})-(\ref{psa2}) admits the two Riemann invariants $c$ and $w=\ln(u)+g(c)$
where $g$ satisfies

\begin{equation}\label{gprime}
 g'(c)    = -\frac{h'(c)}{H(c)}.
\end{equation} 
Indeed, for smooth solutions:
 \begin{eqnarray}  \label{sysentropy}
\left\{ 
  \begin{array} {ccccc}
    \partial_x c &+& \lambda\, \partial_t c   & = & 0, \\
     \partial_x w  & &               & = & 0.
 \end{array}
 \right.
\end{eqnarray}
In the sequel we make use of the  Riemann invariant 
$$W = u\,G(c) = e^w, $$
where
$$ G(c)    =  \exp(g(c)). $$
There are two families of smooth  entropies for the PSA system
(\ref{psa1})-(\ref{psa2}):
  $u\,\psi(c)$ and $\phi(u\,G(c))$ where $\phi$ and $\psi$ are any smooth real functions, and the
corresponding entropy fluxes satisfy 
$$Q'(c)  = h'(c)\, \psi(c) +  H(c)\, \psi'(c).$$ 

 The first family is degenerate convex (in variables  $(u,m)$) provided $\psi''\geq 0$. 
 So, we seek after weak entropy solutions which satisfy, for such functions  $\psi$
 $$ \partial_x \left(u\,\psi(c) \right) + \partial_t Q(c) \leq 0$$ in 
 the distribution sense.
\\
 The second family is not always convex. Nevertheless, if $\pm G'' \geq 0$ then $u\,\psi = \pm
u\,G(c)$ is a degenerate convex entropy,
with entropy flux $ Q\equiv 0$, contained in the family of entropies
$u\,\psi(c)$. So, if $G''$ keeps a constant sign on $[0,1]$, $(c,u)$ has to 
satisfy, in the distribution sense: 
  \begin{equation} \label{IRdecay}
 \pm \frac{\partial}{\partial x} \left(u\,G(c) \right)   \leq  0,
      \quad \mbox{ if } \pm G'' \geq 0 \mbox{ on $[0,1]$.}
\end{equation} 
Notice that the entropies $u\,\psi(c)$ and the entropy $u\,G(c)$ are linear with
respect to the velocity $u$.


\subsection{The Riemann Problem}\label{TRP}

\begin{figure}[H]
\centering
\includegraphics[scale=0.5]{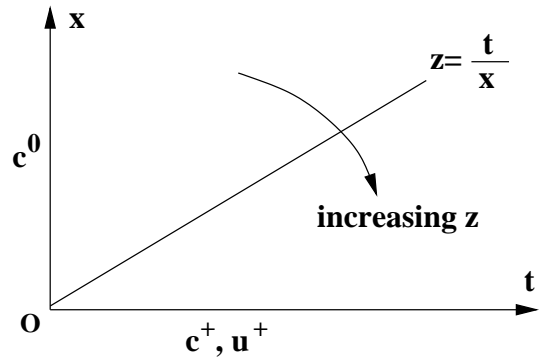}
\caption{data for the boundary Riemann problem \label{drp}}
\end{figure}

For the PSA system it is more convenient to begin by the resolution of the boundary Riemann problem
which is characteristic. The complete Riemann problem is solved at the end of this subsection.
The boundary Riemann problem has the form:
\begin{eqnarray}\label{bPR}
\left\{   \begin{array}{ccc}  \vspace{2mm}
\dx u +\dt h(c)&=&0,\\
\dx  (u c)+\dt I(c)  & = & 0,
\end{array}\right.\label{sysbPR}\\
c(0,x)=c^0 \in [0,1],  \quad x > 0, &\quad
&
\left\{
\begin{array}{ccl}
c(t,0) &=&c^+   \in [0,1],\\
u(t,0)&=&u^+ > 0,
\end{array}
\right.
t>0.\label{databPR}
\end{eqnarray}
We are classically looking for a selfsimilar
solution, i.e.: $\ds c(t,x)=C(z)$,  $u(t,x)=U(z)$ with $z=\ds\frac{t}{x} >0$ (see Fig. \ref{drp}).

In the domain $t>0,\, x>0$, the boundary Riemann problem is solved with $\lambda-$waves since
$\lambda $ is the only positive eigenvalue of the system. 
Let us recall the three following results obtained in \cite{BGJ07}, where $H$ and $f$ are given by
(\ref{eqH}, \ref{eqf}).

\begin{proposition}[ $\lambda-$rarefaction waves]~\\
Any smooth non-constant self-similar solution $(C(z), U(z))$ of (\ref{bPR}) in an open domain\\
$\Omega=\{0\leq\alpha<z<\beta\}$   where $f''(C(z))$ does not vanish, satisfies:
$$
\frac{d C}{ dz}  = \frac{H(C)}{z\,f''(C) },\quad U(z) = \frac{H(C)}{z}.
$$
In particular, $\ds \frac{d C}{ dz}$ has the same sign as $ f''(C)$.\\
Assume for instance  that $0\leq a<c^0<c^+<b\leq 1$ and $f''>0$ in $]a,b[$. Then the only smooth
self-similar solution of (\ref{bPR}) is such that :
\begin{equation}
\left \{\begin{array}{cccr}
C(z)&=&c^0 ,&   0  <  z  <  z_0,\\
\ds\frac{d C}{ dz} & = &
 \ds \frac{H(C)}{z\,f''(C) },&\;   z_0  <   z  < z_+, \\
C(z) & =&   c_+, &   z_+  <  z,
\end{array}  \right.
\end{equation}
where
$ z^+=\ds\frac{H(c^+)}{u^+} $, $z^0=z^+\ds\,e^{-\Phi(c^+)}$ with $\Phi(c)=\ds\int_{c^0}^c
\frac{f''(\xi)}{H(\xi)}\,d\xi$. Moreover $u^0=\ds\frac{H(c^0)}{z^0}$  and $U$ is given by:
\begin{equation}\label{Uz}
\left \{  \begin{array}{cccr}
U(z)& = &   u_0 ,&   0  <   z  <   z_0 ,\\
U(z)& =& \ds\frac{H(C(z))}{z}, & \;  z_0   <   z   < z_+, \\
U(z)& =&   u_+ &   z_+  <   z.
 \end{array}  \right.
\end{equation}
\end{proposition}

\begin{proposition}[$\lambda-$shock waves]\label{SW}
If $(c^0,c^+)$ satisfies the following admissibility condition equivalent to the Liu
entropy-condition (\cite{L76}):
$$\hbox{for all } c \hbox{ between } c^0 \hbox{ and } c^+, \quad \frac{f(c^+)-f(c^0)}{c^+ - c^0}\leq
\frac{f(c)-f(c^0)}{c - c^0},$$
 then the Riemann problem (\ref{bPR}) is
solved by a shock wave defined as:
\begin{equation}
C(z)=\left\{\begin{array}{ccl}
c^0 &\hbox{ if } &  0 < z < s, \\
c^+ & \hbox{ if }&  s<z
\end{array}\right.,
\qquad
U(z)=\left \{\begin{array}{ccl}
u^0 &\hbox{ if }& 0< z < s,\\
u^+ &\hbox{ if } & s < z,
\end{array}  \right.
\end{equation}
 where  $u^0$ and  the speed $s$ of the shock are obtained through 
$$\frac{[f]}{u^0\,[c]} + \frac{1+h^0}{u^0}=s = \frac{[f]}{u^+\,[c]} + \frac{1+h^+}{u^+},$$
where $[c]=c^+ - c^0$, $[f]=f^+ - f^0 = f(c^+) - f(c^0)$, $h^+=h(c^+)$, $h^0=h(c^0)$.
\end{proposition}

\begin{proposition}[$\lambda-$contact discontinuity ]\label{lCD}
Two states $U^0$ and $U^+$ are connected by a $\lambda-$contact discontinuity  with $c^0\neq c^+$
if and only if  $f$ is affine between $c^0$ and $c^+$.
\end{proposition}

In this paper, there is no $\lambda-$contact discontinuity since $f$ is convex: see assumption
$(H_3)$ below.

It appears from these results that we can build a weak entropy solution of the  boundary Riemann 
problem (\ref{bPR}) in a very simple way (see \cite{BGJ07}), similar to the scalar case
with flux $f$, for any data. In particular, if $f''$ has a constant sign, the boundary Riemann
problem is always solved by a simple wave.

\medskip
\noindent
Now, for PSA system \eqref{bPR}, we solve the  Riemann problem  with the following
initial data:
\begin{eqnarray}
\left\{
\begin{array}{ccl}
c(t,0) &=&c^-   \in [0,1],\\
u(t,0)&=&u^- > 0,
\end{array}
\right.
t<0,
&
\left\{
\begin{array}{ccl}
c(t,0) &=&c^+   \in [0,1],\\
u(t,0)&=&u^+ > 0,
\end{array}
\right.
t>0.\label{dataPR}
\end{eqnarray}
A $0-$wave appears on the line $\{t=0\}$.
\begin{proposition}[$0-$contact discontinuity ]\label{0CD}
Two distinct states $U^-$ and $U^0$ are connected by a $0-$contact discontinuity if and only if
$c^-=c^0$ (with of course $u^-\neq u^0$). 
\end{proposition}

In conclusion the solution of the Rieman problem for $x>0$  and a  convex function $f$ is  
\begin{itemize}
 \item $(c,u)= (c^-,u^-)$  for $ t<0$,
 \item a $0-$contact discontinuity  for $t=0$,
 \item a $\lambda-$ wave for $t>0$,
\end{itemize}
see Fig. \ref{Riemann} below.
In practice, since $c^0= c^-$,  we first solve the boundary Riemann Problem \eqref{bPR},
\eqref{databPR}. Then $u^0$ is well defined and the $0-$contact discontinuity
 is automatically solved.

\begin{figure}[H]
\centering
\includegraphics[scale=0.5]{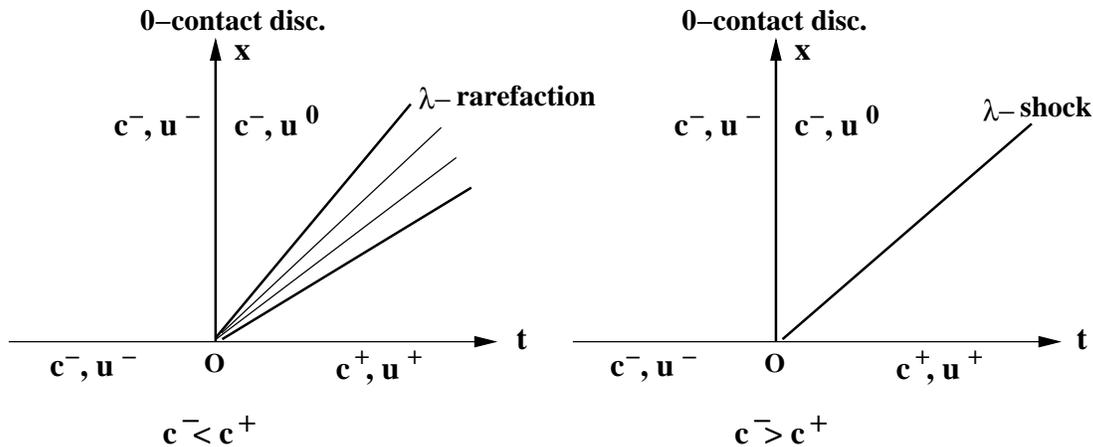}
\caption{solution of the Riemann problem when $f''>0$. \label{Riemann}}
\end{figure}

\begin{figure}[H]
\centering
\includegraphics[scale=0.5]{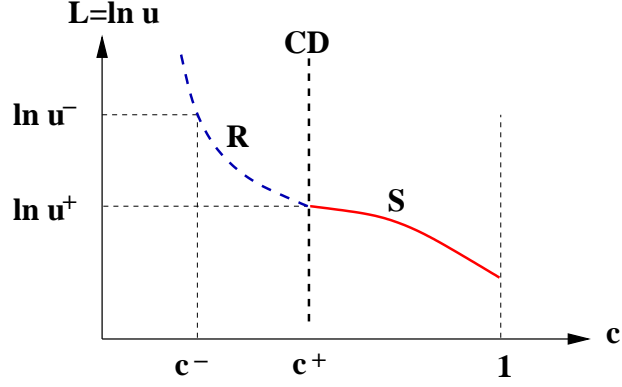}
\caption{structure of shock curves (S), rarefactions (R) and contact discontinuities (CD) in the
case $f''>0>h'$. \label{courbes}}
\end{figure}

We give also (Fig. \ref{courbes}) the structure of shock and rarefaction curves when the function
$f$ is convex, as above,  and $h' < 0$: as stated in \cite{BGJ4}, these curves are monotonic in that
case. Moreover the couple of variables $(L=ln u, c)$ is well adapted to the study of these curves
because $L+g(c)$ is a Riemann invariant.

\subsection{Main assumptions}\label{mainass}

PSA system has the strong property to admit a positive Riemann invariant $W=u\,G(c) > 0 $, which is
also an entropy. This Riemann invariant  plays a key-role in the blow up mechanism: only provided we
ensure that $W$ is non decreasing with respect to $x$ we can expect, in view of (\ref{maxp}), the
velocity $u$ to increase and blow up. Indeed,  if $W$ is decreasing with respect to $x$ we have
$0\leq u\,G(c)\leq u_b(t)\,G(c_b(t))$, thus $\ds 0\leq u\leq \Vert
u_b\Vert_\infty\,\sup_{[0,1]}G/\inf_{[0,1]}G$ and $u$ is bounded. Then our  first assumption is:
\begin{eqnarray*} 
 G'' & <& 0  \qquad \qquad \qquad  \qquad \qquad \qquad  (H1).
\end{eqnarray*}
$(H1)$ means, as seen above (\ref{IRdecay}),  that $-W$ is an admissible degenerate convex entropy
for weak entropy
solutions (with a zero entropy-flux) i.e. $\partial_x W\geq 0 $. 

In contrary, notice that in some cases, for instance an inert gas and an active gas with a
Langmuir isotherm  as in \cite{BGJ06}, we have $G'' > 0$, then $W$ is bounded and there is no blow
up for the velocity in $L^\infty$.

The next hypothesis is useful to the construction of our blow up example in a classical
hyperbolic framework of Lax with monotonic shock and rarefaction curves as in \cite{BGJ4}: 
\begin{eqnarray*} 
 h' & \neq & 0, \mbox{ everywhere }  
 \qquad \qquad \qquad \qquad \qquad \qquad  (H2).
\end{eqnarray*}
Assumption $(H2)$ means that one gas is more active than the other. Mathematically, it implies the
monotonicity of $g$ and $G$ because $g'=-h'/H$. 
Thanks to $(H2)$ rarefaction  curves, which are the level line of $W$, are  monotonic curves in
the plane $(c,u)$.
\\
We  also assume that the eigenvalue $\lambda=H(c)/u$ is genuinely nonlinear. Our next assumption is
then:
\begin{eqnarray*} 
 f'' & \neq & 0, \mbox{ everywhere }  
    \qquad \qquad  \qquad \qquad \qquad \qquad (H3).
\end{eqnarray*}

For Temple systems (see the next section for a definition) it is well known that there is no blow up
in $L^\infty$, then we have to
investigate if the PSA system is in this class or not: this is done in the next section. This will
lead us to state this supplementary assumption to avoid the Temple class:
\begin{eqnarray*} 
 \mbox{ {\it PSA system (\ref{psa1})-(\ref{psa2}) is not a Temple system} } \qquad (H4).
\end{eqnarray*}
These four hypothesis are not natural, in the sense that they are not generally satisfied, but some
applications are concerned, as the case of an inert gas associated with ammonia or water vapor for
instance. 
\section{Temple class and PSA system} \label{secTemple}
%
It is well known (\cite{BrG00, Bi01}) that blow up cannot occur for Temple systems with $L^\infty$
data. The aim of this section is thus to investigate the conditions for the  PSA system
(\ref{psa1})-(\ref{psa2}) to be or not to be in the Temple class. Of course, in the sequel, we must 
avoid the Temple class to construct our blow up solution.

About Temple systems, on may consult    \cite{T83-1,T83-2,S87-1,H94,S96}. Let us recall
the definition of the Temple class (\cite{T83-2}).

\begin{definition}[Temple class]
 A system of conservation laws in one space dimension is in the Temple class if it is strictly
hyperbolic and satisfies the following properties:
\begin{enumerate}
\item there exists a system of coordinates consisting of Riemann invariants,
\item shock and rarefaction curves coincide.
\end{enumerate}
\end{definition}

For  PSA system we have the following criterion:

\begin{lemma} \label{temple} 
PSA system (\ref{psa1})-(\ref{psa2}) is a Temple system if and only if $\partial_x W=0$ for the
entropy solution of any Riemann problem (\ref{bPR}), \eqref{dataPR}.
\end{lemma}

\pr 
as seen above, PSA system is strictly hyperbolic and satisfies the first condition of the
definition. There is a linearly degenerate eigenvalue (zero) and the rarefaction curves are only
given by $\partial_x W=0$ (see
(\ref{sysentropy}) or \cite{BGJ07}). 

Notice that the Riemann problem for PSA system has a unique piecewise smooth entropy solution. If
PSA system is in the Temple class, the second condition in the definition ensures that $\partial_x
W=0$ everywhere. The converse is immediate.
\cqfd

Let us highlight a sufficient condition for the PSA system to be in the Temple class:

\begin{proposition} \label{Gsecnul}
If $G''=0$ then PSA system (\ref{psa1})-(\ref{psa2}) is a Temple system.
\end{proposition}

\pr according to (\ref{IRdecay}), if $G''\equiv 0$  we have $\partial_x (uG(c))
= 0$ and from Lemma \ref{temple}, PSA system is a Temple system. \cqfd

We precise now some simple examples where $G''=0$.

\begin{proposition} \label{prop32}
If the two  isotherms are linear then $G''=0$. 
\end{proposition}
\pr
we have $G=\exp(g), \, G'=g' \exp(g), \, G''=(g''+g'^2 )\exp(g)$. With $g'=-h'/H$, 
$H=1+q'_1-ch'(c)$ and $h=q_1+q_2$ we get
$$G''=\displaystyle \frac{\exp(g)}{H^2}(-q_1''-q_2''-q_2''q_1'+q_2'q_1'').$$
For two linear isotherms we have $q_1''=q_2''=0$ and then $G''=0$.
\cqfd

More generally we have:

\begin{lemma}
Assume (H2), then $G''=0$ if and only if there exist two real constants $\alpha,\,\beta$ such that
$$\alpha\,q_1 + (\alpha-1)\,q_2+c+\beta=0. $$
\end{lemma}

\pr as seen above,  $G''=0$ if and only if $q_1''+q_2''+q_2''q_1'-q_2'q_1''=0$. Setting $h_1=q'_1,
\, h_2=q'_2$ this condition writes $(h_1+h_2)'=(h_1+h_2)'h_2-h_2'(h_1+h_2)$ or also
$$\frac{(h_1+h_2)'}{(h_1+h_2)^2}=\frac{(h_1+h_2)'h_2-h'_2(h_1+h_2)}{(h_1+h_2)^2}.$$
Thus  $G''=0$ if and only if there exists a constant $\alpha$ such that
$$\displaystyle \frac{-1}{h_1+h_2}=-\frac{h_2}{h_1+h_2}+\alpha,$$
i.e. $-1=-q'_2+\alpha(q'_1+q'_2)$, which gives immediately the condition. 
\cqfd

In the particular case where a gas is inert, say the first one ($q_1=0$), we can give a very simple
criterion:

\begin{proposition}
Assume $q_1=0$. Then PSA system (\ref{psa1})-(\ref{psa2}) is a Temple system if and only if the
``active'' isotherm is linear ($q_2''=0$).
\end{proposition}

\pr  if $q_2''=0$ and $q_1=0$ we have $G''=0$ thanks to Proposition \ref{prop32}. Then the system is
in the Temple class thanks to Proposition \ref{Gsecnul}.

Reciprocally, assume that the system is in the Temple class: following Lemma \ref{temple}, 
$\partial_x
(uG(c))=0$ holds for the entropy solution of any Riemann problem (\ref{bPR}). Let  $c_-,\,
c_+,\, u_+$ be the data of the Riemann problem: the unknown is $u_-$ and using the
Rankine-Hugoniot conditions we have (see \cite{BGJ07}):

\begin{equation}\label{srp}
u_-=\mathcal{R}(c_-,\, c_+,\, u_+)=\left\lbrace
\begin{array}{rl}
 u_+\,\gamma(c_-,c_+)& \hbox{for a shock}\\
 u_+\,G(c_+)/G(c_-)& \hbox{for a rarefaction,}
\end{array} \right.
\end{equation} 

where
$$\gamma(c_-,c_+)=\frac{[c+q_1(c)]-c_+[h]}{[c+q_1(c)]-c_-[h]}=\displaystyle\frac{\alpha+h_-}{
\alpha+h_+}$$
with $\alpha=\displaystyle \frac{[f]}{[c]}+1$, $[c]=c^+ - c^-$, $[f]=f^+ - f^- = f(c^+) - f(c^-)$
(and so on).

Since the system is in the Temple class, shock and rarefaction curves coincide, thus we have:
$\displaystyle\frac{\alpha+h_-}{\alpha+h_+}=\frac{u_-}{u_+}=\frac{G_+}{G_-}$, then $-\alpha\,
[G]=[hG]$. We can assume $g(0)=0$, then choosing  $c_-=0$  and setting $c=c_+$, we
get $[G]=G-1$, $[hG]=h\,G-K$ with $K=h(0)>0$ and $\alpha=1-h(c)$. Thus, the equality
$-\alpha\,[G]=[hG]$, reads $G=1+K-h$.

Differentiating this relation we get  $g'\,G=-h'$, then $G=\displaystyle \frac{-h'}{g'}=H$.
Differentiating $G=H$ we get $g'=\displaystyle \frac{H'}{H}$, but in view of (\ref{gprime}) we have
$g'=\displaystyle \frac{-h'}{H}$ thus $h'=-H'=h'+ch''$ and finally $h''=0$ i.e. $q_2''=0$.
\cqfd

\begin{remark}\label{ex}
For an inert gas (for instance the first one: $q_1=0$) and an active gas with strictly convex
isotherm ($q_2'' >0$), PSA system is not in the Temple class. For example, it is the case if the
active gas is the ammonia or the  water vapor. For other examples, see \cite{BGJ07,BGJ4}.
\end{remark}
It is now clear that (H4) can be satisfied in some physically relevant cases.

\section{The blow up solution}\label{secBlowup}

In this section we present and prove our main result: the existence of an entropy blow up solution.
It is organized as follows: first we state and comment this result, next we give a basic step
allowing to increase the velocity through two non interacting  waves, a shock and a rarefaction,
then we repeat this process to build a piecewise Lipschitz entropy solution with growing
velocity. Lastly, a sequence of such solutions leads to a blow up (in a bounded space-time domain).

\subsection{Main result}

Our main result needs four assumptions explained in subsection \ref{mainass}.

\begin{theorem}[Blow up for  velocity] \label{ThBUV}
Assume (H1), (H2), (H3) and (H4). Then for any  $T>0, \, X_\infty>0$, there exists a set of  initial
and boundary data (\ref{psa0}) and a corresponding weak entropy solution on $[0,T] \times
[0,X_\infty[$ of PSA system (\ref{psa1})-(\ref{psa2}) such that 
$$\mid\mid u\mid\mid_{L^\infty((0,T)\times (0,X_\infty))} =+\infty.$$
\end{theorem}

Actually the solution  built to prove this theorem has special features.

\begin{itemize}
\item The velocity  only blows up at the boundary $\{t=0\}$ when  $x \rightarrow X_\infty$, with
arbitrary small data.
\item The system is strictly hyperbolic, but at ( $t=0, \, x=X_{\infty}$) it becomes degenerate
hyperbolic because $\lambda =\ds \frac{H(c)}{u}$ tends to $0$ when $u$ tends to $+\infty$. Then the
boundary becomes twice characteristic.
\item
    $\forall X\in ]0, X_\infty[,\,  u,\,c \in L^\infty(0,T; BV(0,X))
\cap L^\infty(0,X; BV(0,T))$. The concentration $c$ remains bounded while $u$ blows up.

\item
Let $\Omega$ be a neighborhood of the critical point $(t=0,x=X_\infty)$ 
 such that $\Omega \subset [0,T] \times [0,X_{\infty}]$.
 Outside $\Omega$, we can prove that $(u,c)$ 
has a piecewise smooth structure, so the blow up occurs only at the boundary.
Indeed, there is an accumulation of wave-interactions near  $(t=0,x=X_\infty)$. 
\item
To build such solutions, we have necessarily to choose the boundary concentration $c_0(x)=c(0,x)$ 
 not in $ BV(0,X_\infty)$. Else there is no blow up, see  \cite{BGJ06,BGJ07}.

\end{itemize}

\subsection{A mechanism to increase velocity}\label{secConstruction}

In this section, we present the main ingredient for the construction of the blow up solution. It
consists in the resolution of two consecutive boundary Riemann problems which leads to increase the velocity
without increasing the concentrations.

\begin{figure}[H]
\centering
\includegraphics[scale=0.4]{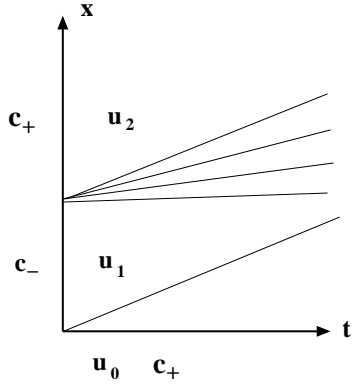}
\caption{the two boundary Riemann problems \label{trp}}
\end{figure}

For the first boundary Riemann problem the data are $(c_-,\, c_+,\, u_0)$ chosen in such a way that the
solution is
a shock wave and we set (see (\ref{srp})) $u_1=\mathcal{R}(c_-,\, c_+,\, u_0)$. For the second
problem the data are $(c_+,\, c_-,\, u_1)$, the solution is necessarily a rarefaction wave and we
set  $u_2=\mathcal{R}(c_+,\, c_-,\, u_1)$.

The sign of $G''$ is crucial  but not those of  $f''$ and $h'$ in assumptions $(H3)$, $(H2)$.
Indeed, exchanging the labels of the two gas changes signs of $f''$ and $h'$.  \textbf{In the
sequel, we fix these signs:}

$$ f''< 0,  \quad h' < 0 \quad \mbox{\rm  everywhere} $$ 
but it is not restrictive and all variants work. 

With these assumptions we must take $c_-<c_+$ (see Propositions \ref{SW} and \cite{BGJ06}).

We introduce now  the \textit{amplification coefficient} $\bR$ defined by 
$$u_2=\bR u_0. $$
We are going to show that this coefficient only depends upon $(c_-,c_+)$.

\begin{lemma}
The amplification  coefficient  satisfies $\bR=\bR(c_-,c_+)$.
\end{lemma}
\pr according to \eqref{srp}, for the shock the Rankine-Hugoniot condition gives $${u_1}=
u_0\,\gamma(c_-,c_+)$$
and for the rarefaction, we have $\partial_x (u\,G(c))=0$ then $u_2\,G(c_+)=u_1\,G(c_-)$. Finally we
have $$u_2=\frac{G(c_-)}{G(c_+)}\,\gamma(c_-,c_+)\,u_0= \bR(c_-,c_+)\,u_0.$$
\cqfd

The following lemma justifies the label ``amplification''.

\begin{lemma}\label{Rgeq1}
We assume $(H1)$, then the amplification  coefficient satisfies  $\bR \geq 1$.
\end{lemma}
\pr since $G'' < 0$ then  $\partial_x (u\,G(c)) \geq 0$ (entropy condition (\ref{IRdecay})).
Through the shock wave we have $u_1 > u_0$,  $c_- < c_+$ and $u_0\, G(c_+) \leq u_2\, G(c_+)$.
Finally, $u_2=\bR u_0 \geq u_0$ then $\bR \geq 1$.\cqfd

Notice that for the shock curves we have $c_- < c_+$ and $u$ is not monotonic along the
process because $u_0\leq  u_2\leq u_1$. Actually we need  a strict amplification i.e.
$\bR>1$. We show that it is true for almost every choice $(c_-,\,c_+)$ with $c_-<c_+$, provided the
isotherms are analytic as we will see below. We proceed in two steps.

\begin{lemma}
$\bR(c_-,c_+) \equiv 1$  if and only if the system is in the Temple class.
\end{lemma}
\pr assume that $\bR(c_-,c_+)=1$ then $u_2=u_0$. This is true for any $(c_-,\,c_+)$ with $c_-<c_+$
and
thus the shock curve connecting $(c_+,u_0)$ to $(c_-,u_1)$ and the  rarefaction curve connecting
$(c_-,u_1)$ to $(c_+,u_2=u_0)$ are the same, but traveled in opposite directions. Finally shock
and rarefaction curves coincide i.e the system is in the Temple class.

Reciprocally, if the system is a Temple system we have  $\partial_x (u\,G(c))=0$, then
$$u_0\,G(c_+) =u_1\,G(c_-) =u_2\,G(c_+), $$
finally $u_0=u_2$ and $\bR(c_-,c_+)=1$.
\cqfd

\begin{lemma}
We assume $(H1)$, $(H3)$ and $(H4)$.\\ If the function $\bR$ is analytic then $\bR(c_-,c_+) > 1$ 
almost everywhere in the domain $c_- < c_+$.
\end{lemma}
\pr
thanks to $(H1)$ and lemma \ref{Rgeq1} we have $\bR\geq 1$. Since $\bR$ is analytic then either
$\bR(c_-,c_+)=1$ always, either $\bR(c_-,c_+) > 1$ almost everywhere. With the assumption $(H4)$,
PSA system is not a Temple system then, thanks to the previous lemma, $\bR > 1$ almost everywhere. 
\cqfd

Notice that the  commonly used isotherms are analytic, then $\bR$ also. With Remark \ref{ex} we
have some relevant examples where the situation $\bR>1$ occurs.

\subsection{A solution with growing velocity} \label{ssPLS}

Following the mechanism introduced in the preceding subsection we impose constant  initial data and
piecewise constant concentration at the boundary such that Riemann problems at the boundary are
alternatively solved by a shock or a rarefaction.

Let $N > 1$ be a fixed integer, $0=x_0 < x_1 < \cdots < x_{2N-1} < X=x_{2N}$, 
$0< \underline{c} <   \overline{c}<1$ such that $\bR( \underline{c}, \overline{c}) > 1$ and $u_0 >
0$. 
We solve PSA system (\ref{psa1})-(\ref{psa2}) with the following data for $0 < t < T$, $0 < x < X$,
$k=0,1,\cdots,N-1$:
\begin{eqnarray} \label{datapc} 
 \left\{
 \begin{array}{ccc}
  c(t,0)   & = &     \overline{c},  \\
  u(t,0)   & = &     u_0, \\
  c(0,x)  & = & 
    \ds \left \{ 
             \begin{array}{ccc}
               \underline{c}  & if & x_{2k} < x < x_{2k+1}, \\
               \overline{c}  & if & x_{2k+1} < x < x_{2k+2}. 
             \end{array}
  \right.
\end{array}
\right.
\end{eqnarray}

With these data we prove existence and uniqueness of a piecewise Lipschitz entropy solution:

\begin{theorem}
We assume $(H1)$, $(H2)$, $(H3)$.
With da\-ta (\ref{datapc}), there exists a unique weak entropy solution in the class of piecewise
Lipschitz functions. Furthermore, this solution has exactly $N$ shock curves on $[0,T] \times
[0,X]$.
\end{theorem}

The uniqueness is an important feature to extend the previous entropy solution until the blow up.
 
\begin{figure}[H]
\centering
\includegraphics[scale=0.5]{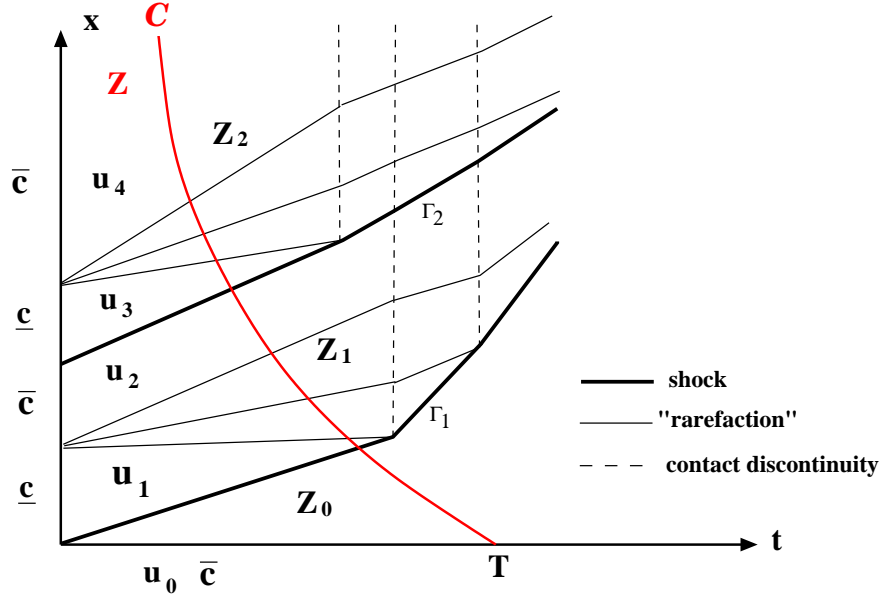}
\caption{first interactions and free domain Z \label{Z}}
\end{figure}

\pr 
to construct the solution we use the Front Tracking Algorithm as in \cite{BGJ4}
(see $\cite{Br00,D00}$  for general references). There exists a domain $Z$ (see Fig. \ref{Z}), a
neighborhood of the boundary $\{t=0\}$ where boundary Riemann problems do not interact. In this domain $Z$ we
have an explicit
solution. Let us denote by $u_k$ the value of $u(x,0^+)$ when $ x_{k-1} < x < x_{k}$ for a given $k
>0$. 

Since $f''<0$ and  $h'<0$ we have $N$ shocks emerging from $((x_{2k},t=0))_{k=0}^{N-1}$ and $N$
rarefactions from $((x_{2k+1},t=0))_{k=0}^{N-1}$. Furthermore $ u_{2k} < u_{2k+1}$, $u_{2k}<
u_{2k+2}=
\bR\, u_{2k}$ and 

\begin{equation}\label{Rk}
u_{2k}= \bR^k\, u_0.
\end{equation}

Take $\delta >0$ to discretize the rarefactions as in $\cite{BGJ4}$ and let shocks, rarefactions and
contact discontinuities interact to obtain an approximate entropy solution on $[0,T]\times[0,X]$.

As in a scalar conservation laws with piecewise constant data, no new shock can appear but shocks
can disappear. It is a consequence of wave-interactions studied in \cite{BGJ4}. We recall the
following results concerning interactions:

 \begin{itemize}
\item  if a rarefaction interacts with a shock then we have a contact discontinuity and a
rarefaction or a shock. 

If the shock and the rarefaction have the same strength (i.e. the same
absolute variation of the concentration through the simple wave) we have only a contact
discontinuity;
\item 
if two shocks interact we obtain a contact discontinuity and a stronger shock;
\item 
if a shock interacts with a contact  discontinuity, we obtain a contact discontinuity and a shock
with the same strength;
\item  
if a rarefaction interacts with a contact discontinuity, we obtain a contact discontinuity and a
rarefaction with the same strength.    
\end{itemize}

So, we can follow each shock, more precisely:

\begin{itemize}
\item 
if a shock disappears after an interaction with a stronger rarefaction,
by convention we follow the characteristic speed associated with the eigenvalue $\lambda=H(c)/u$,
and the strength of the shock is $0$;
\item 
if two shocks interact, they become a single shock. But, we consider this shock as two shocks to
define without ambiguity the parametrization of the $N$ shocks issuing from 
$((x_{2k},t=0))_{k=1}^N$ by $N$ functions $s_k^\delta :[0,X] \rightarrow [0,T]$. 
\end{itemize}
The entropy inequality  $\partial_x u\,G(c) \geq 0$  yields a  lower bound $\inf u >0$ (see
\cite{BGJ06}): indeed we have $u\,G(c)\geq u_b(t)\,G(c_b(t))$, $\ds 
G(c_b)\geq \inf_{[0,1]}G>0$, $\ds 0< G(c)\leq \sup_{[0,1]}G$ and $ u_b(t)\geq\alpha>0$. Then we get
an upper bound for $\lambda=H(c)/u$: $0 <\lambda < \bar{\lambda}=\ds\frac{\sup_{[0,1]} H}{\inf
u}$. Thus the speeds of the  shock curves $\ds \frac{d s_k^\delta }{dx}$ are bounded by
$\bar{\lambda} $.

Thus, we have $N$ shock curves uniformly Lipschitz. When  $\delta$ goes to $0$, Ascoli's Theorem
gives $N$ Lipschitz curves $(s_k)_{k=1}^{N}$, up to a subsequence. We denote by
$U^\delta=(c^\delta,u^\delta)$ the approximate solution issued from the Front Tracking Algorithm
with $\delta>0$ fixed. Extracting again a subsequence $(\delta_n)$, the sequence $U_n=U^{\delta_n}$
converges towards an entropy solution $U=(u,c)$ as in \cite{BGJ06,BGJ07,BGJ4}.
Let us denote by $\bS$ the set of shocks emerging from the boundary.  
$$  \bS= \ds \bigcup_{k=1}^N \Gamma_k  \qquad \hbox{ where } \Gamma_k=\{ (x,s_k(x)),\, x \in
[0,X]\}$$ 
and $\Gamma_k^\delta$, $S^\delta$ with a similar definition  for the approximate solution
$U^\delta$.   Let us show that $c$ has no singularity except on $\bS$.

\medskip
\underline{Outside $\bS$ the concentration $c$ is necessarily continuous:}
notice that $c$ has $BV$ regularity (see \cite{BGJ07,BGJ4}), then $c$ is discontinuous  if and only
if $c$ has shocks. Furthermore $(c,u)$ is  an entropy solution hence  $\partial_t c >0 $
through the shocks. So  it suffices to prove that $\partial_t c \leq 0$ outside $\bS$ to
conclude. 

For this purpose, take any open set $\Omega $ such that  $S \cap \Omega = \emptyset$.
Using the uniform limit of the shock curves, we have no shock for $U_n$ in $\Omega$ and  $n$
sufficiently large: $ n > n_1$.  So, in $\Omega$,  any jump of $c$ satisfies
 $ [c^{\delta_n}] \leq 0$ since it is in a rarefaction domain $\Omega$ for $n > n_1$.
That is to say $ \partial_t c^{\delta_n} \leq 0 $ in $\Omega$ for $n > n_1$.
Passing to the limit we get $\partial_t c  \leq 0$.
Then, $c$ is continuous in $\Omega$. This is true for any open set outside $\bS$,
so $c$ is continuous except on $\bS$. 

Since $c$ is continuous outside a finite number of shocks, we can obtain  $c$ and $u$ by the method
of characteristics, as in \cite{BGJ4}. More precisely, it was proved  in \cite{BGJ4} that any
entropy solution with a continuous concentration is unique and can be constructed by the method of
characteristics. Since $\lambda$ is bounded, the characteristics are Lipschitz continuous, so $c$ is
Lipschitz too in the complementary of $\bS$, which prove the piecewise Lipschitz structure of the
concentration for the entropy solution. Notice that $u$ is less regular than $c$ because there are
contact discontinuities emerging at each interaction between two waves.

\medskip
\underline{Uniqueness:}
the piecewise Lipschitz  structure allows us to prove uniqueness. It is a simple extension of a
theorem proved in \cite{BGJ06} for piecewise smooth solutions. For completeness, we give the key
points of this proof. Under the first shock, namely in the domain $Z_0$ (see Fig. \ref{Z}), $(u,c)$
is unique and determined by the method of characteristics. In our case, $(u,c)$ is constant in
$Z_0$, but we do not use this fact to give a general argument and to use it for the other domains.  

Over the first shock and under the second one, namely domain $Z_1$, the characteristics are well
defined if we know $u$ on the upper side of $\Gamma_1$ since $\partial_x u\,G(c)=0$ on $Z_1$. Indeed
we have the following equations in variables $(u,c)$ for some functions $\mathcal{F}$ and
$\mathcal{S}$:
$$
\begin{array}{lccr}
\partial_x W  & = & 0, &    \mbox{ in }  Z_0 \cup Z_1, \\    
  \partial_t c + W \,\phi(c)\, \partial_x c  & = & 0,  
       &    \mbox{ in }  Z_0 \cup Z_1, \\
   \partial_t \chi(t,s,x)  & = & W \phi(c),  
        &  \mbox{ in }  Z_0 \cup Z_1, \\
   W_1-W_0 & = &\mathcal{F}(c_0,c_1),  &    \mbox{ on }   \Gamma_1, \\
   \ds\frac{ds_1}{dt} &=&\mathcal{S}(W_0,c_0,c_1), 
 &   \mbox{ on }  \Gamma_1
\end{array}
$$
with $W(t,0), c(t,0),c(0,x) $ are given functions and  $\phi(c)=1/(H(c)G(c))$. 

The trick in \cite{BGJ06} is to work in variable $t$. In the same way, $t$ is an admissible
parameter for the shock curves $\Gamma_k$: $\ds\frac{ds_k}{dt}$ is well defined since we have BV
estimates for $\ln u$, thus $\lambda$ avoids zero. Moreover $W$ has jumps with respect to $t$
(contact discontinuities for $u$) and the product $\ds\frac{1}{W\,\phi(c)}\,\partial_t c$ may not be
defined.

The first two equations are equations for the Riemann invariants $W$ and $c$. The third equation is
the equation of characteristics associated to the eigenvalue $\lambda$ where $\chi(t,t,x)=x$. The
fourth equation is the Rankine-Hugoniot equation, where $c_i$ and $W_i$ are the traces of $c$ in
$Z_i$ on
 $\Gamma_1$. This relation  gives us $W$ in $Z_1$ if we know $\Gamma_1$, $c_0$ and $c_1$. Indeed
$\Gamma_1$ is an unknown curve, but the last equation determines it with an ordinary differential
equation. 

At $t=0$, we know $W_0,c_0,c_1$ so we know
 $\Gamma_1$ and $W_1$ locally near the corner for all $t < \tau$, and then for all $x $. We use a
local uniqueness theorem as the Cauchy-Lipschitz Theorem. Notice that when we know $W_1$ for $ t <
\tau$, the characteristics are also defined for all $x $ thanks to $\partial_x W=0$. 
  So, we have a uniqueness result in $Z_0\cup \Gamma_1 \cup Z_2$ for $t < \tau$ and  we extend the
uniqueness result until $t=T$. 
 Next, we use the same technique for the second shock curve and so on.
\\
Now, as in \cite{BGJ06},  $U=(c,u)$ is unique. So, the whole sequence $U^\delta$ tends to $U$.
\cqfd

\subsection{Blow up of the velocity} \label{ssBU}

Now, we use data (\ref{datapc}) with $N=+\infty$ i.e. $(x_k)_{k \in \N}$ is an increasing sequence
such that $\ds \lim_{k \rightarrow + \infty} x_k = X_\infty$.
Then, the concentration remains bounded at the boundary but the $BV$ norm of $c$ and the $L^\infty$
norm of $u$ blow up at the boundary $t=0$.
Near $(t=0, x = X_\infty)$, there is an accumulation of oscillations for the concentration and an
accumulation of interactions between shocks and rarefactions.\\
Let  $Z$ be the subset of  $[0,T] \times [0,X_\infty[$ defined in the previous subsection: it is a
neighborhood of the vertical segment $\{0\} \times [0,X_\infty[$.
Indeed, the first interaction takes place at $(t^1,x^1)$: since a contact discontinuity propagates
vertically, we have first to cut the set $\{t > t^1, x > x^1\}$. We do the same work, for all
first interactions of the $2N$ Riemann problems issuing from the boundary.

Then $ Z =\{0 < x <X_\infty, \quad 0 < t < z(x)\}$ where the function $ z: [0,X_\infty] \rightarrow
[0,T]$ is piecewise constant on $[0,X]$ for any $X < X_\infty$, non increasing, positive before
$X_\infty$,  $ z(0)=T$,  $z(X_\infty)=0.$

On $Z$, we exactly know  the solution and 
  $$\lim_{(t,x) \in Z, \; x \rightarrow X_\infty} u(t,x) = + \infty.$$

{\it \bf Proof of Theorem \ref{ThBUV}:}
let $0 < X< X_{\infty}$.
By previous construction, we get a unique piecewise smooth entropy solution $U^X$ on $[0,T] \times
[0,X]$. If $0 < X < Y < X_{\infty}$, we get $U^Y$ and by uniqueness, $U^X = U^Y$ on $[0,T] \times
[0,X]$. So, as for an ordinary differential equation, we can consider the unique maximal solution on
$[0,T] \times [0,X_{\infty}[$.
Before interaction between the solutions of the  Riemann problems issuing from the boundary, we
know explicitly $U$ in $Z$. From geometric increasing of $u(0,x)$ when $x \rightarrow X_{\infty}$,
we get
a blow up for $u$ at $t=0$, $x =X_{\infty}$.

Furthermore, the characteristic slope $1/\lambda$ goes to infinity near $\{0\}\times\{X_\infty\}$
but outside a suitable ``triangular neighborhood'' of this point we get a determination domain where
the solution has BV regularity and there is no blow up. This determination domain has the form
$$\mathcal T=\{ (t,x)\,;\,0<t<T,\ 0<x<\min(X+\frac{t}{\overline{\lambda}},X_{\infty}\}$$
for any $X$ in $]0,X_{\infty}[$.
\cqfd

\section{Conclusion}
 To conclude this paper, let us give some comments  about this topic.  
 \begin{itemize}
  \item  It is a new  physical  example of  hyperbolic partial differential  equations  exhibiting a
blow up. However, this chemical  model becomes  not physically relevant for too large velocity, $u
>> 1$.

 \item It follows from our blow up example that general $L^\infty$ theory for PSA system  only with
$c_0, c_b, \ln u_b \in L^\infty$ is impossible. But, since our example does not blow up at finite
distance of the boundary, a natural conjecture is: ``if $\ln u$ remains bounded at the boundary
$t=0$, there is no blow up''.

\item We also conjecture that our solution can be defined for $x\geq X_\infty$ as an entropy
locally piecewise smooth solution outside the boundary, with infinite velocity on $\{0\}\times
[X_\infty,+\infty[ $.
 \end{itemize}
  
 Finally we give  some comparisons  with other blow up solutions  for $2\times 2$ systems. The
reader can  also consult \cite{SY05} and the memoir  \cite{S07}.
 \begin{itemize}
  \item  There is a very simple and impressive example in  \cite{Y03short,Y03}, with two Burgers
equations linearly coupled on two boundaries. This solution  blows up in finite time $T$ for all
$x$. In contrary, for our example, there is only one boundary and the solution only blows up on
this  boundary.
  \item As in  \cite{JY04}, it is an example exhibiting a blow up of
amplitude  where at the same time the system looses strict hyperbolicity. This is a special feature
of $2\times 2$ systems, since for larger systems, at least $3\times 3$ systems, this is well known
that blow up can occur without loosing the strict hyperbolicity.

  \end{itemize}
 

\begin{thebibliography}{10}

\bibitem{Bi01}
S.~Bianchini.
\newblock Stability of ${L}^\infty$ solutions for hyperbolic systems with
  coinciding shocks and rarefactions.
\newblock {\em SIAM J. Math. Anal.}, 33(4):959--981, 2001.

\bibitem{B92}
C.~Bourdarias.
\newblock Sur un syst\`{e}me d'edp mod\'{e}lisant un processus d'adsorption
  isotherme d'un m\'{e}lange gazeux. (french) [on a system of p.d.e. modelling
  heatless adsorption of a gaseous mixture].
\newblock {\em M2AN}, 26(7):867--892, 1992.

\bibitem{B98}
C.~Bourdarias.
\newblock Approximation of the solution to a system modeling heatless
  adsorption of gases.
\newblock {\em SIAM J. Numer. Anal.}, 35(1):13--30, 1998.

\bibitem{BGJ06}
C.~Bourdarias, M.~Gisclon, and S.~Junca.
\newblock Some mathematical results on a system of transport equations with an
  algebraic constraint describing fixed-bed adsorption of gases.
\newblock {\em J. Math. Anal. Appl.}, 313(2):551--571, 2006.

\bibitem{BGJ07}
C.~Bourdarias, M.~Gisclon, and S.~Junca.
\newblock Existence of weak entropy solutions for gas chromatography system
  with one or two actives species and non convex isotherms.
\newblock {\em Commun. Math. Sci.}, 5(1):67--84, 2007.

\bibitem{BGJ08}
C.~Bourdarias, M.~Gisclon, and S.~Junca.
\newblock Hyperbolic models in gas-solid chromatography.
\newblock {\em Bol. Soc. Esp. Mat. Apl.}, 43:29--57, 2008.

\bibitem{BGJ4}
C.~Bourdarias, M.~Gisclon, and S.~Junca.
\newblock Strong stability with respect to weak limit for a hyperbolic system
  arising from gas chromatography. arXiv:0907.1733v2. 
\newblock {\em Submitted in Methods and Applications of Analysis}.

\bibitem{Br00}
A.~Bressan.
\newblock {\em Hyperbolic Systems of Conservation Laws. The One-Dimensional
  Cauchy Pro\-blem}, volume~20 of {\em Oxford Lecture Series in Mathematics and
  its Applications}.
\newblock Oxford University Press, 2008.

\bibitem{BrG00}
A.~Bressan and P.~Goatin.
\newblock Stability of ${L}^\infty$ solutions of {T}emple class systems.
\newblock {\em Differential Integral Equation}, 13(10-12):1503--1528, 2000.

\bibitem{CGM03}
C.~Cheverry, O.~Gu\`es, and G.~M\'etivier.
\newblock Oscillations fortes sur un champ lin\'eairement d\'eg\'en\'er\'e.
  (french) [strong oscillations on a linearly degenerate field].
\newblock {\em Ann. Sci. Ecole Norm. Sup.}, (4) 36(5):691--745, 2003.

\bibitem{CGM04}
C.~Cheverry, O.~Gu\`es, and G.~M\'etivier.
\newblock Large-amplitude high-frequency waves for quasilinear hyperbolic
  systems.
\newblock {\em Adv. Differential Equations}, 9(7--8):829--890, 2004.

\bibitem{CG99}
A.~Corli and O.~Gu\`es.
\newblock Local existence of stratified solutions to systems of balance laws.
\newblock In {\em Ann. Univ. Ferrara Sez.}, volume VII, pages 109--119,
  Workshop on Partial Differential Equations, Ferrara, 1999.

\bibitem{CG01}
A.~Corli and O.~Gu\`es.
\newblock Stratified solutions for systems of conservation laws.
\newblock {\em Trans. Amer. Math. Soc.}, 353(6):2459--2486, 2001.

\bibitem{D00}
C.~Dafermos.
\newblock {\em Hyperbolic Conservation Laws in Continuum physics}.
\newblock Springer, Heidelberg, 2000.

\bibitem{H94}
A.~Heibig.
\newblock Existence and uniqueness of solutions for some hyperbolic systems of
  conservation laws.
\newblock {\em Arch. Ration. Mech. Anal.}, 126:79--101, 1994.

\bibitem{JY04}
H.-K. Jenssen and R.~Young.
\newblock Gradient driven and singular flux blowup of smooth solutions to
  hyperbolic systems of conservation laws.
\newblock {\em J. Hyperbolic Differ. Equ.}, 1(4):627--641, 2004.

\bibitem{DCRBT88}
M.~Douglas Levan, C.A. Costa, A.E. Rodrigues, A.~Bossy, and D.~Tondeur.
\newblock Fixed--bed adsorption of gases: Effect of velocity variations on
  transition types.
\newblock {\em AIChE Journal}, 34(6):996--1005, 1988.

\bibitem{L76}
T.P. Liu.
\newblock The entropy condition and the admissibility of shocks.
\newblock {\em J. Math. Anal. Appl.}, 53:78--88, 1976.

\bibitem{M04}
A.~Museux.
\newblock Stratified weak solutions of the 1-d lagrangian euler equations are
  viscosity solutions.
\newblock {\em Adv. Differential Equations}, 9(11--12):1395--1436, 2004.

\bibitem{RAA70}
H.~K. Rhee, R.~Aris, and N.R. Amundson.
\newblock On the theory of multicomponent chromatography.
\newblock {\em Philos. Trans. R. Soc. Lond.}, (267):419--455, 1970.

\bibitem{RAA86}
H.~K. Rhee, R.~Aris, and N.R. Amundson.
\newblock {\em First-Order Partial Differential Equations}, volume~I.
\newblock Prentice-Hall, 1986.

\bibitem{RSVG88}
P.~Rouchon, M.~Sghoener, P.~Valentin, and G.~Guiochon.
\newblock Numerical simulation of band propagation in nonlinear chromatography.
\newblock In {\em Chromatographic Science Series}, volume~46, New York, 1988.
  Eli Grushka, Marcel Dekker Inc.

\bibitem{Ru84}
P.M. Ruthwen.
\newblock {\em Principles of adsorption and adsorption processes}.
\newblock Wiley Interscience, 1984.

\bibitem{S87-1}
D.~Serre.
\newblock Solutions \`a variations born\'ees pour certains syst\`emes
  hyperboliques de lois de conservation.
\newblock {\em J. Differential Equations}, 68:137--168, 1987.

\bibitem{S96}
D.~Serre.
\newblock {\em Syst\`emes de lois de conservation I}.
\newblock Diderot, Paris, 1996.

\bibitem{S07}
M.~Sever.
\newblock {\em Distribution solutions of nonlinear Systems of conservation
  Laws}.
\newblock Memoirs of the AMS, 2007.

\bibitem{SY05}
W.~Szeliga and R.~Young.
\newblock Blowup with small bv data in hyperbolic conservation laws.
\newblock {\em Arch. Ration. Mech. Anal.}, 179:31--54, 2005.

\bibitem{T83-1}
B.~Temple.
\newblock Systems of conservation laws with coinciding shock and rarefaction
  curves.
\newblock In J.~Smoller, editor, {\em Nonlinear Partial Differential
  Equations}, volume~17, pages 143--151, Providence, RI, 1983. Amer. Math. Soc.

\bibitem{T83-2}
B.~Temple.
\newblock Systems of conservation laws with invariant manifolds.
\newblock {\em Trans. Amer. Math. Soc.}, 280:781--795, 1983.

\bibitem{Y99}
R.~Young.
\newblock Exact solutions to degenerate conservation laws.
\newblock {\em SIAM J. Math. Anal.}, 1999.

\bibitem{Y03short}
R.~Young.
\newblock Blowup in hyperbolic conservation laws.
\newblock {\em Contemporary Mathematics}, 1(2):269--292, 2003.

\bibitem{Y03}
R.~Young.
\newblock Blowup of solutions and boundary instabilities in nonlinear
  hyperbolic equations.
\newblock {\em Commun. Math. Sci.}, 1(2):269--292, 2003.

\end{thebibliography}

\end{document}